\documentclass[a4paper, 12pt, leqno]{article}

\usepackage{amsmath,amsfonts}
\usepackage[latin1]{inputenc}
\usepackage[active]{srcltx}
\begin{document}
\bibliographystyle{plain}
\newtheorem{theo}{Theorem}[section]
\newtheorem{lemme}[theo]{Lemma}
\newtheorem{cor}[theo]{Corollary}
\newtheorem{defi}[theo]{Definition}
\newtheorem{prop}[theo]{Proposition}
\newtheorem{problem}[theo]{Problem}
\newtheorem{remarque}[theo]{Remark}
\newcommand{\beq}{\begin{eqnarray}}
\newcommand{\enq}{\end{eqnarray}}
\newcommand{\be}{\begin{eqnarray*}}
\newcommand{\en}{\end{eqnarray*}}
\newcommand{\Td}{\mathbb T^d}
\newcommand{\Rd}{\mathbb R^n}
\newcommand{\R}{\mathbb R}
\newcommand{\N}{\mathbb N}
\newcommand{\Sn}{\mathbb S}
\newcommand{\Snn}{\mathbb S^{n-1}}
\newcommand{\Zd}{\mathbb Z^d}
\newcommand{\Linf}{L^{\infty}}
\newcommand{\dt}{\partial_t}
\newcommand{\Dt}{\frac{d}{dt}}
\newcommand{\Dtt}{\frac{d^2}{dt^2}}
\newcommand{\demi}{\frac{1}{2}}
\newcommand{\vf}{\varphi}
\newcommand{\epu}{_{\epsilon}}
\newcommand{\ep}{^{\epsilon}}
\newcommand{\bfi}{{\mathbf \Phi}}
\newcommand{\bpsi}{{\mathbf \Psi}}
\newcommand{\bx}{{\mathbf x}}
\newcommand{\ds}{\displaystyle}
\newcommand{\mT}{\text{c-exp}}
\newcommand{\bT}{\mathbf{T}}
\newcommand{\cexp}{\exp^c}

\let\cal=\mathcal
\title{Regularity of optimal maps on the sphere:
\\  The quadratic cost and the reflector antenna}
\author{Gr\'egoire Loeper\footnotemark[1]}
%\date{}
\maketitle
%\begin{center}
%Gr\'egoire LOEPER\footnotemark[1]
%\end{center}
\footnotetext[1]{Institut Camille Jordan, Universit\'e Claude Bernard Lyon 1}

\begin{abstract}
Building on the results of Ma, Trudinger and Wang \cite{MTW}, and of the author \cite{L5},  we study two problems of optimal transportation on the sphere: the first corresponds to the cost function $d^2(x,y)$, where $d(\cdot,\cdot)$ is the Riemannian distance of the round sphere;
the second  corresponds  to the cost function $-\log|x-y|$, it is known as the reflector antenna problem.
We show that in both cases, the {\em cost-sectional curvature} is uniformly positive, and establish the geometrical properties so that  the results of \cite{L5} and \cite{MTW} can apply: global smooth solutions exist for arbitrary smooth positive data and optimal maps are H\"older continuous under weak assumptions on the data. 
\end{abstract}

\newpage
%\tableofcontents

\section{Introduction}

\subsection{Monge-Kantorovitch problem on a Riemannian manifold}
Let  $M$ be a  topological space, 
let $c: M\times M\to \R \cup \{+\infty\}$ be a cost function and $\mu_0, \mu_1$ be probability measures on $M$. 
In the optimal transportation problem, one looks for a map $\bT: M \to M$ that minimizes the functional
\beq\label{min}
{\cal I}(\bT) = \int_M c(x, \bT(x))d\mu_0(x),
\enq
under the constraint that $\bT$  pushes forward $\mu_0$ onto $\mu_1$, (hereafter $\bT_{\#}\mu_0 = \mu_1$), i.e.
\be
\forall B\subset M \;\text{Borel},\,\, \mu_1(B) = \mu_0(\bT^{-1}(B)).
\en

This problem has been first studied by Monge \cite{Monge}, with $M= \Rd$, for the cost $c = |x-y|$ (the Euclidean distance). Beyond Monge's original problem, Brenier studied the case of the
quadratic cost $c=|x - y|^2$, and pointed out its close connection with important nonlinear PDEs (Monge-Ampere, Euler etc...). For the quadratic cost, when $\mu_0$ is absolutely continuous with respect to the Lebesgue measure, he proved the existence and uniqueness of an optimal map T. This maps has a convex potential (i.e. $T=\nabla \phi$ with $\phi$ convex) and is shown to be the only map with convex potential that  pushes forward $\mu_0$ onto $\mu_1$.
After Brenier's result, the theory of optimal transportation has been extended to general cost functions. The existence of optimal maps is granted under very generic conditions on the cost function, and the way to obtain it is achieved through a general procedure,  known as  Kantorovitch duality: Optimal maps are obtained by solving the dual Monge-Kantorovitch problem, whose unknown are potential functions. For $\phi$ a lower semi-continuous function on $M$, we define its {\it c-transform} as
\be
\phi^c(x) &=& \sup_{y\in M}\{-c(x,y) - \phi(x)\}.
\en
A potential $\phi$ is c-convex if it is the c-transform of some other $\psi : M \to \R$. In that case, the equality $\phi=\phi^{cc}$ holds. (Notice that the quadratic cost is equivalent to the cost $-x\cdot y$, for which the c-transform is nothing but the Legendre-Fenchel transform; hence, $[-x\cdot y]$-convex functions are convex functions.)
Under suitable assumptions, and following for example \cite{CaGenCost}, the minimizers in (\ref{min}) are related to c-convex potentials as follows: for an optimal $\bT_{opt}$ in (\ref{min}), there exists a c-convex potential $\phi$ such that 
\beq\label{defGphi}
\text{for a.e.} \,x \in  M, \bT_{opt}(x)=G_\phi(x):=\{y \in M, \phi(x) + \phi^c(y)   = -c(x,y)\}.
\enq 
Conversely, if $\bT: M\to M$ can be expressed under the form (\ref{defGphi}) for some c-convex $\phi$, for $\mu_0$ a probability measure on $M$ and $\mu_1$ its push-forward by $\bT$, then $\bT$ is the optimal map between $\mu_0$ and $\mu_1$.
Of course, it is not clear a-priori that (\ref{defGphi}) defines a map, as $G_\phi(x)$ is a set. However, under a suitable assumption on the cost (assumption {\bf A1} below), the set $G_\phi(x)$ will be reduced to a single point for Lebesgue almost every $x$. Note also that when $M$ is compact (which we will assume throughout the remainder of the paper), the set $G_\phi(x)$ is never empty when $\phi$ is c-convex.

Brenier's result was generalized in a natural way to Riemannian manifolds by McCann: Let $M$ be a  manifold,
with Riemannian metric $g$, compact and without boundary, with distance function $d(\cdot,\cdot)$. For $u,v \in T_x(M)$, $(u,v)_g(x)$ (or in short $(u,v)_g$) denotes the scalar product on ${T}_x(M)$ with respect to the metric $g$, $|v|_g^2=(v,v)_g$.
From the results of \cite{Mc2}, in the case where $c = d^2/2$, the optimal map can be expressed as a so-called  gradient map, i.e.
\be
G_\phi(x) = \exp_x(\nabla_g\phi(x)),
\en
where $\nabla_g$ denotes the gradient with respect to  the Riemannian metric $g$ on $M$ (and from now, we omit the subscript $g$), and $\phi$ is some c-convex potential. For a general cost, 
one needs first to introduce the c-exponential map $\mT_x(\cdot)$, defined as the inverse of $y \to -\nabla_xc(x,y)$ (again see assumption {\bf A1} below).
Optimal maps will then be given  by 
\be
\bT_{opt}(x)=G_\phi(x):= \mT_x(\nabla\phi(x)).
\en
(This definition holds in the a.e. sense and is consistent with the definition (\ref{defGphi}).)
In all cases, for a smooth potential $\phi$ such that $G_{\phi\,\#}\mu_0 = \mu_1$, the conservation of mass is expressed in local coordinates by the Monge-Amp\`ere type equation
\beq
\label{ma}
\det(D^2\phi + D^2_{xx}c(x,G_\phi))= \frac{\rho_0}{\rho_1(G_\phi)}|\det D^2_{x,y}c(x,G_\phi)|,
\enq
where $\rho_0, \rho_1$ are densities of $\mu_0, \mu_1$ with respect to the Lebesgue measure.

A complete description of the optimal transportation problem can be found in \cite{Vi}, otherwise the introduction of \cite{L5} encloses the necessary material for the present paper. We also mention the second book by Villani on optimal transport \cite{Vi2}, which presents some of the results enclosed in this paper.

In the present work, we will address the problem of regularity of solutions of (\ref{ma}) (or equivalently of minimizers of (\ref{min})) in the particular geometrical setting of the constant curvature sphere  of $\Rd$, that we will denote $\Snn$. We will consider two cases: the quadratic cost $c(x,y)=\demi d^2(x,y)$ with $d$ the Riemannian distance, and the reflector antenna case $c(x,y)=-\log|x-y|$ (see below). We will show that in both cases, the regularity results obtained in \cite{L5} and \cite{MTW} hold. As noticed at the end of the paper, our result can be easily generalized to the case $c(x,y) = f(d(x,y))$, assuming some conditions on $f$ (see Theorem \ref{main_general}).

The problem of regularity of optimal maps is related to the problem of regularity of solutions of the associated elliptic Monge-Amp\'ere type equation. This fully non-linear elliptic partial differential equation has received considerable attention over the last decades. Up to recently, the only case where regularity results had been obtained were for the Monge-Amp\`ere equation
\be
&&\det D^2 u= f,\\
&& u \text{ convex}, 
\en
with works by Caffarelli \cite{Ca1, Ca0, Ca2, Ca4, Ca3, Ca5}, Urbas \cite{U}, and Delano\"e \cite{De}.
This form of the Monge-Amp\`ere equation is associated to the quadratic cost. It was only recently with the works of Ma, Trudinger and Wang \cite{MTW, TW} and subsequent results by the author \cite{L5} that regularity results for generic costs were obtained.
 Our goal in the present paper is to continue this study in the Riemannian setting, in the particular case of the round sphere.

\paragraph{The reflector antenna problem}
Consider a closed hypersurface $\Sigma$ of $\Rd$ parametrized by a so-called "height" function $h$ from $\Snn$ to $\R$, i.e. $\Sigma = \{xh(x), x\in \Snn\}$. Given $\Omega, \Omega'$ two domains of $\Snn$, and two probability measures $\mu_0, \mu_1$ on $\Omega, \Omega'$,
the reflector antenna problem is then to find $h$, under the constraint that the antenna reflects the incoming  intensity $\rho_0$ into the outgoing  of intensity $\rho_1$.
If a ray with direction $x$ is reflected into a ray of direction $T(x)$, the conservation of energy turns into a Monge-Amp\`ere type equation for $h$ (see \cite{CaGuHu}):
\be
\frac{\det\left(\nabla_{ij}h + (u-\eta)e_{ij}\right)}{\eta^{n-1}\det e_{ij}}=\frac{\rho_0(x)}{\rho_1(T(x))}
\en
with $$\eta = \frac{|\nabla h|^2 + h^2}{2h},$$
$\nabla_{ij}h$ the second covariant derivative of $h$, and $e_{ij}$ the Riemannian metric of $\Snn$. 

Existence, regularity and uniqueness of weak and strong solutions have been addressed by several authors, among them Wang \cite{Wa3}, Guan and Wang \cite{GuWa}, Glimm and Oliker \cite{GlimOlik1}, and Oliker \cite{Olik}.  The connection with optimal transport theory was established by Wang in \cite{Wa4}, who showed that the reflector antenna problem is equivalent to an optimal transport problem on the sphere with cost equal to $-\log|x-y|$.
 
\paragraph{Organization of the paper}
In the next section, we will expose our results; then we will give a reminder of previous results and give some definitions and notations (some of them can be needed to understand the results). The rest of the paper is dedicated to the proofs of the results.

\section{Results}

This work addresses the issue of regularity of optimal maps. We deal only with global solutions (the measures and transport maps are defined on the whole of $\Sn^{n-1}$).
When the data are positive and $C^2$ (resp. $C^\infty$) smooth, we show that the optimal potential is $C^3$ (resp. $C^\infty$).
When the target measure is bounded by below, and the source measure satisfies {\bf B(n-1)}, see (\ref{mini}) (resp. {\bf A(n-1,p)} for some $p>n-1$, see (\ref{Lp})) the optimal potential is $C^1$ (resp. $C^{1,\alpha}$ for $\alpha(n,p)$). 
We also give an original and self-contained proof of the connectedness of the contact set $G_\phi(x)$ (defined in (\ref{defGphi})) in those particular cases. Our results are identical for the quadratic cost $c(x,y)=\demi d^2(x,y)$ and the reflector antenna $c(x,y) = -\log|x-y|$. 

For the reflector antenna problem, classical smooth solutions of (\ref{ma})  had been obtained in \cite{Wa3, Wa4, GuWa}. Global $C^1$ regularity for weak solutions of (\ref{ma}) had also been obtained independently in \cite{CaGuHu}, under the assumption that both the source and target measures have densities  bounded away from 0 and infinity. Here we relax their assumptions, allowing the source measure to vanish, and requiring an integrability condition that does not even imply absolute continuity with respect to the Lebesgue measure. Moreover the H\"older exponent in the $C^{1,\alpha}$ result is explicit.

Our results are not a simple corollary of the results of \cite{L5} once the positivity of the cost-sectional curvature is asserted. Indeed, the basic assumption {\bf A0} is not satisfied  on $\Sn^{n-1} \times \Sn^{n-1}$, due to singularities of the cost functions (cut-locus), hence the equation (\ref{ma}) itself could become singular.
We have to show first that the graph of optimal transport map lies in a subdomain of $\Sn^{n-1} \times \Sn^{n-1}$ which is uniformly far from the cut locus. This is done by improving a result established in \cite{del-loep}, and also adapting it to the antenna case (see the Propositions \ref{sphereok} and \ref{awayantenna}). 

Then we can localize the problem and reduce it locally to an Euclidean problem, and thus use  the $C^{1,\alpha}$ estimates of \cite{L5} to show partial regularity  under assumption (\ref{Lp}) or (\ref{mini}).

For classical regularity, we  employ the method of continuity: combining the results of \cite{del2} with the crucial a-priori estimate established in \cite{MTW}, we obtain classical smooth solutions. 

We choose to present in details the proof of the results for the quadratic case, from which the antenna case follows easily, once the key ingredients are verified.

\paragraph{Acknowledgments}
I wish to thank Alessio Figalli and Cédric Villani for precious remarks and fruitful discussions.   I also thank Philippe Delano\"e with whom we started to think about the problem of regularity for optimal transportation on the sphere, and Robert McCann who first raised to me the issue of the connectedness of the contact set, in 2003. 
I gratefully acknowledge the support of a French Australian exchange grant PHC FAST EGIDE No.12739WA. I wish also to thank Neil Trudinger, Xu-Jia Wang and the Center for Mathematics and its Applications at University of Canberra for their hospitality, as part of this work was accomplished there.

\subsection{The quadratic cost}

In this part, we consider  $\Sn^{n-1}$ the unit sphere of $\Rd$ equipped with the round metric $g$, and Riemannian distance $d$, and we let $c(x,y)=\demi d^2(x,y)$. We first have the following essential remark:
\begin{prop}\label{essential}
For all $x \in \Snn$, the set $\{p \in \partial_x c(x,y), y\in \Snn\}$ is the closed ball   $\bar B(0,\pi)$. Letting $\hat x$ be the antipodal point of $x$, the set $\{-\nabla_xc(x,y), y\in \Snn\setminus \hat x\}$
is the open ball $B(0,\pi)$. Hence for all $y_0, y_1 \in \Sn^{n-1} \setminus \hat x$, the c-segment $[y_0,y_1]_x$ (see definition \ref{def-c-seg}) is well defined.
\end{prop}

The next result shows that {\bf As} is satisfied outside of the cut-locus.
\begin{theo}\label{Assphere}
Let $\rm{antidiag}$ be the set $\{(x,\hat x), x\in \Sn^{n-1}\}$ with $\hat x$ the antipodal point of $x$.
Then the cost-sectional curvature is uniformly positive on $\Sn^{n-1}\times \Sn^{n-1}\setminus \rm{antidiag}$.
\end{theo}

Then we show that the contact set $G_\phi(x)$ is always connected (and even c-convex). The definitions of $\partial\phi(x), \partial^c\phi(x)$ are reminded below in Definitions \ref{sub-dif}, \ref{c-sub-dif}.
\begin{theo}\label{main_sphere_1}
Let $\phi$ be a c-convex potential on $\Sn^{n-1}$. For all $x\in \Sn^{n-1}$, let the set $G_\phi(x)$ be defined as in (\ref{defGphi}). Then 
\begin{enumerate}
\item $G_\phi(x)$ is c-convex with respect to $x$, and $\partial\phi(x) = \partial^c\phi(x)$.
\item When $\hat x \in G_\phi(x)$, $G_\phi(x)$ is equal to the whole sphere $\Sn^{n-1}$.  
\end{enumerate}
\end{theo}

Finally, we conclude with our regularity result:
\begin{theo}\label{main_sphere}
Let $\mu_0,\mu_1$ be two probability measures  on $\Sn^{n-1}$ and $\phi$ be a c-convex potential. Assume that $G_{\phi\#} \mu_0  = \mu_1$. Assume that $\mu_1 \geq m{\rm dVol}$, for some $m>0$. Then
\begin{enumerate}
\item If $\mu_0$ satisfies {\bf B(n-1)} (see (\ref{mini})), then $\phi\in C^1(\Sn^{n-1})$, and as in Theorem \ref{main}, the modulus of continuity of $\nabla\phi$ depends on $f$ in (\ref{mini}).
\item If $\mu_0$ satisfies {\bf A(n-1,p)} (see (\ref{Lp})) for some $p>n-1$, then $\phi \in C^{1, \beta}(\Sn^{n-1})$ with $\beta=\beta(n-1,p)$ as in Theorem \ref{main}.
\item If $\mu_0, \mu_1$ have positive $C^{1,1}$ (resp. $C^\infty$)  densities with respect to  the Lebesgue measure, then $\phi \in C^{3,\alpha}(\Sn^{n-1})$ for every $\alpha \in [0,1[$ (resp. $\phi \in C^\infty(\Sn^{n-1})$.)
\end{enumerate}
\end{theo}

\subsection{The reflector antenna}
The results presented above adapt with almost no modification to the reflector antenna.
Before that, we remark that by changing $y$ into $-y$, it is equivalent to study the cost $c(x,y) = -\log|x+y|$. Then, as in the previous case, the set of singular points of the cost is equal to $\rm{antidiag}$ (i.e. the set of antipodal pair of points).  Moreover, it is straightforward to check that whenever $c(x,y)$ satisfies {\bf A0}, {\bf A1}, {\bf A2}, {\bf AS}, then $c(x,-y)$ also satisfies those assumptions. Then, we prove the
\begin{theo}\label{mainantenna} 
Theorem \ref{main_sphere_1} and \ref{main_sphere} hold for $c(x,y) = -\log|x - y|$. 
\end{theo} 
(For point 2 of Theorem \ref{main_sphere_1}, notice that one must not consider $\hat x$ the antipodal point of $x$, but $x$ itself.)
Before entering into the proofs of our results, we present a review of previous results and concepts from the work of Ma, Trudinger and Wang in \cite{MTW} and the author in \cite{L5}. 

\section{Cost-sectional curvature : from geometry of contact sets to elliptic regularity}
This paragraph is a short reminder of the results established in \cite{L5}, \cite{TW}, \cite{MTW}. 
In a nutshell, the things to remember are the following: when the contact set defined in (\ref{defGphi}) is unconditionally connected, the associated Monge-Amp\`ere equation behaves well (i.e. smooth data lead to smooth solutions). Moreover, there is a tensor whose non-negativity is equivalent to the connectedness of the contact set. This tensor is the Ma-Trudinger-Wang tensor, or {\it cost-sectional curvature} tensor. When the cost is the squared Riemannian distance, this tensor contains informations  about the Riemannian curvature of the underlying manifold.   

We adapt the results of \cite{L5}, \cite{TW}, \cite{MTW} to the Riemannian case by considering considering $M$ a $n$-dimensional manifold, and $D$ a general domain of $M\times M$ (which for simplicity we assume compact), instead of restricting to tensorial domains of $\Rd \times \Rd$, i.e. domains of the form $\Omega\times\Omega'$ for $\Omega, \Omega'$ domains of $\Rd$. 
We denote $\pi_1, \pi_2$ the usual canonical projections. For any $x\in \pi_1(D)$, we denote by $D_x$ the set $D\cap \pi_1^{-1}(x)$. We proceed similarly for the $y$ variable. 
Let us introduce the following conditions:
\begin{itemize}
\item[{\bf A0}] The cost function $c$ belongs to $C^4(D)$.
\item[{\bf A1}] For all $x \in \pi_x(D)$, the map $y\to -\nabla_xc(x,y)$ is injective on $D_x$. 
\item[{\bf A2}] The cost function $c$ satisfies $\det D^2_{x,y}c \neq 0$ for all $(x,y)$ in $D$.
\end{itemize}
Assumption {\bf A1} allows to introduce the following definition:
\begin{defi}\label{defiT} 
Under assumption {\bf A1}, for $x\in \pi_1(D)$ we define the c-exponential map at $x$, which we denote by $\mT_x$, the map such that 
\be
\forall (x,y)\in D, \mT_x(-\nabla_xc(x,y))=y.
\en
\end{defi}
The c-exponential map coincides with the Riemannian exponential map when $c = d^2/2$. In the general case,  optimal maps are defined by
\be
\bT_{opt}(x)=G_\phi(x)=\mT_x(\nabla \phi(x))
\en
for some c-convex potential $\phi$.

Noticing that for all $y\in G_\phi(x)$, $\phi(\cdot) + c(\cdot, y)$ has a global minimum at $x$, we introduce /  recall  the following definitions:
\begin{defi}[subdifferential]\label{sub-dif}
For $\phi$ a semi-convex function, the subdifferential of $\phi$ at $x$, that we denote $\partial\phi(x)$, is the set 
\be
\partial\phi(x)=\Big\{p\in T_x(M): \forall y \in M, \phi(y) \geq \phi(x) + (p\cdot \exp_x^{-1}(y))+ o(d(x,y))\Big\}.
\en
\end{defi}
The subdifferential is always a convex set, and is always non empty for a semi-convex function. 
\begin{defi}[c-subdifferential]\label{c-sub-dif}
If $\phi$ is $c$-convex, the c-sub-differential of $\phi$ at $x$, that we denote $\partial^c\phi(x)$, is the set
\be
\partial^c\phi(x)=\Big\{ -\nabla_x c(x,y), y\in G_\phi(x)\Big\}.
\en
The inclusion $\emptyset \neq \partial^c\phi(x)\subset \partial\phi(x)$ always holds.
\end{defi}
We make here the following important remark: $\mT_x(\partial^c\phi(x))$ is set of all $y$ such that $-\phi(\cdot) -c(\cdot,y)$ reaches a global maximum at $x$,  while $\mT_x(\partial\phi(x))$ is the set of all $y$ such that $-\phi(\cdot) -c(\cdot,y)$ has a critical point at $x$.

We remind the  definition of c-convexity (see \cite{MTW}): 
\begin{defi}[c-convex sets]\label{def-ccvx-set}
Let $x\in \pi_1(D)$. A subset $\omega'$ of $\pi_2(D_x)$  is c-convex (resp. uniformly  c-convex) with respect to $x$ if the set $\{-\nabla_xc(x,y), y \in \omega'\}$ is a convex (resp. uniformly convex) set of $T_xM$. 
Whenever $\omega \times \omega' \subset D$,  $\omega'$  is c-convex with respect to  $\omega$  if it is c-convex with respect to every $x \in \omega$.
\end{defi} 
There exits also a notion of segment associated to the cost $c$, which we call c-segment:
\begin{defi}[c-segment] \label{def-c-seg} Let $(x,y_0)$ and $(x,y_1)\in D$. Let $p_i = -\nabla_xc(x,y_i), i=0,1$. Assume that 
$[p_0, p_1] \subset -\nabla_xc(x,D_x)$. Then the c-segment between $y_0$ and $y_1$ with respect to $x$, which we denote $[y_0, y_1]_x$ is defined by
\be
[y_0, y_1]_x=\big\{y\in D_x, -\nabla_xc(x,y) \in [p_0, p_1]\big\}.
\en
Then $y_\theta$ will be such that $-\nabla_xc(x,y_\theta) = \theta p_1 + (1-\theta) p_0$.
\end{defi}

We now remind the definition of the {\em cost-sectional curvature}.
Although this name and the present formulation  comes from \cite{L5}, this tensor has first been introduced by Ma, Trudinger and Wang in \cite{MTW}.
\begin{defi}\label{defiR}
Under assumptions {\bf A0}-{\bf A1}-{\bf A2},  for every $(x_0,y_0) \in D$, one can define  on $T_{x_0} M \times T_{x_0} M$ the real-valued map
\beq\label{defR}
 \  \  \  \  \   \mathfrak{S_c}(x_0,y_0)(\xi,\nu)= D^4_{p_\nu p_\nu x_\xi x_\xi}\Big[(x,p)\to -c(x,\mT_{x_0}(p))\Big]\Big|_{x_0,p_0=-\nabla_x c(x_0,y_0)}.
\enq
When $\xi,\nu$ are unit orthogonal vectors (with respect to the metric $g$ at $x_0$), $\mathfrak{S_c}(x_0,y_0)(\xi,\nu)$ defines  the  cost-sectional curvature from $x_0$ to $y_0$ in directions $(\xi,\nu)$. This definition is intrinsic (i.e. coordinate independent).
\end{defi}
We are now ready to introduce the following condition on the cost function:
\begin{itemize}
\item[{\bf As}] \textbf{(Positive sectional curvature)} The cost-sectional curvature is uniformly positive on $D$.
In other words, there exists $C_0 > 0$ such that for all $(x_0, y_0) \in D$, for all $\xi, \nu \in T_{x_0}M$, $(\xi,\nu)_g=0$,
\end{itemize}
\beq\label{C0}
\mathfrak{S_c}(x_0,y_0)(\xi,\nu) \geq C_0|\nu|^2|\xi|^2.
\enq

Alternatively, condition {\bf Aw} \textbf{(non-negative sectional curvature)} is satisfied whenever (\ref{C0}) is satisfied with $C_0=0$.

In \cite{MTW}, \cite{TW}, \cite{L5} the study was restricted to the case where $M= \Rd$, $D = \Omega \times \Omega'$, for $\Omega, \Omega'$ two domains of $\Rd$. 
Concerning classical regularity, Ma, Trudinger and Wang \cite{MTW} proved the following:
\begin{theo}[\cite{MTW}]\label{classicalmtw}
Let $c$ be a $C^4(\Omega \times\Omega')$ cost function that satisfies assumptions  {\bf A1}, {\bf A2}, {\bf As} on $(\Omega\times\Omega')$, $\Omega'$ being  c-convex with respect to $\Omega$. 
Let $\mu_0 = \rho_0 dx, \mu_1 = \rho_1 dx$ be probability measures  respectively on $\Omega$ and $\Omega'$.
Assume that $\rho_0, \rho_1$ are bounded away from 0, and belong to $C^2(\Omega)$ (resp. $C^2(\Omega')$). 
Let $\phi$ be a c-convex potential on $\Omega$ such that $G_{\phi\,\#}\mu_0=\mu_1$.
Then $\phi \in C^{3,\alpha}(\Omega)$ for every $\alpha \in [0,1)$.
\end{theo}

The key step toward this result is an interior a-priori estimate for second derivatives of solutions of (\ref{ma}); we will rely on this a-priori estimate in the present paper.

In \cite{L5} we proved the following relationship between regularity and a natural geometric property of the c-convex potentials:

\begin{theo}[\cite{L5}]\label{main0}
%Assuming {\bf A0} {\bf A1}, and the c-convexity of $\Omega, \Omega'$ with respect to each other. 
Assume that for all $\mu_0, \mu_1$ smooth positive probability measures, the c-convex potential $\phi$ such that $G_{\phi\#}\mu_0 = \mu_1$ is $C^1$. Then for all $\phi$ c-convex on $\Omega$, for all $x \in \Omega$, $\partial\phi(x) = \partial^c\phi(x)$.
\end{theo}

This result, combined with Theorem \ref{classicalmtw} yields easily the following:

\begin{theo}[\cite{L5}]\label{main1}
Let $c$ be a $C^4(\Omega \times\Omega')$ cost function that satisfies assumptions  {\bf A1}, {\bf A2}, {\bf As} on $(\Omega\times\Omega')$, $\Omega, \Omega'$ being  c(resp c*)-convex with respect to each other. 
Then, for all $\phi$ c-convex on $\Omega$, for all $x\in \Omega$, the set $G_\phi(x)$ is connected.
\end{theo}

In \cite{L5}, it was also proved that  the non-negativity of the cost-sectional curvature (i.e. condition {\bf Aw}) is a necessary condition for regularity: without that condition, one can construct potentials that are not $C^1$ even though the data are smooth and positive. Relying on the results of \cite{TW}, it was then shown that the non-negativity of the cost-sectional curvature was equivalent to the connectedness of the set $G_\phi(x)$ (see (\ref{defGphi})).

\begin{theo}[\cite{L5}, \cite{TW}]
Let $c$ be a $C^4(\Omega \times\Omega')$ cost function that satisfies assumptions  {\bf A1}, {\bf A2} on $(\Omega\times\Omega')$, $\Omega, \Omega'$ being  uniformly c(resp c*)-convex with respect to each other.
Then the following are equivalent:
\begin{enumerate}
\item The cost $c$ satisfies {\bf Aw};
\item For all $\phi$ c-convex, for all $x\in \Omega$, $G_\phi(x)$ is connected;
\item For all $\mu_0, \mu_1$ $C^\infty$-smooth positive probability measures on $\Omega, \Omega'$,  the optimal map between $\mu_0$ and $\mu_1$ is $C^\infty$.
\end{enumerate} 
\end{theo} 

Under positive cost-sectional curvature, we proved also $C^{1,\alpha}$ regularity of optimal potentials for rough data.
Let us first introduce some regularity properties on the measures $\mu_0, \mu_1$.
The first one reads
\begin{itemize}
\item[\bf A(n, p)]For some $p\in ]n, +\infty], C_{\mu_0}>0$,
\end{itemize}
\beq\label{Lp}
\mu_0(B_\epsilon(x)) \leq C_{\mu_0} \epsilon^{n (1-\frac{1}{p})} \text{ for all } \epsilon \geq 0, x \in M.
\enq  
The second condition reads
\begin{itemize} 
\item[\bf B(n)]For some $f:\R^+\to \R^+$ with $\lim_{\epsilon \to 0} f(\epsilon) = 0$, 
\end{itemize}
\beq\label{mini}
\mu_0(B_\epsilon(x)) \leq f(\epsilon) \epsilon^{n (1-\frac{1}{n})} \text{ for all } \epsilon \geq 0, x \in M.
\enq
We remark that $\mu\in L^p$ implies {\bf A(n, p)} with the same $p$, while {\bf B(n)} doest not imply $\mu \in L^1$.
Letting $\cal{ H}^{n-1}$ be the $n-1$ dimensional Hausdorff measure, we also notice that (\ref{mini}) implies that $\mu$ does not give mass to sets $A$ such that $\cal{ H}^{n-1}(A)$ is finite, which is close to the optimal assumption for existence of an optimal map in (\ref{min}) (e.g. $\mu_0$ does not charge  $n-1$-rectifiable sets).
We also denote classically $\Omega_\delta = \{x\in\Omega: d(x,\partial\Omega)>\delta)$. 

\begin{theo}[\cite{L5}]\label{main}
Under the assumptions of Theorem \ref{main1} on $c$, $\Omega$, $\Omega'$,
let $\mu_0, \mu_1$ be probability measures  respectively on $\Omega$ and $\Omega'$.   Let $\phi$ be a c-convex potential on $\Omega$ such that $G_{\phi\,\#}\mu_0=\mu_1$.
Assume that $\mu_1  \geq m\,{\rm dVol}$ on $\Omega'$ for some $m>0$.
\begin{enumerate}
\item  Assume that $\mu_0$ satisfies {\bf A(n,p)}  for some $p>n$.
Let $\alpha = 1-\frac{n}{p},  \  \beta = \frac{\alpha}{4n-2+\alpha}$.
Then for any $\delta>0$ we have
\be
\|\phi \|_{C^{1,\beta}(\Omega_\delta)} \leq {\cal C},
\en 
and ${\cal C}$ depends only 
on $\delta>0$, $C_{\mu_0}$ in  (\ref{Lp}), on $m$,  on the constants  in conditions {\bf A0}, {\bf A1}, {\bf A2}, {\bf As}.

\item If $\mu_0$ satisfies {\bf B(n)},
then $\phi$ belongs to $C^1(\Omega_\delta)$ and the modulus of continuity of $\nabla\phi$ depends also on $f$ in (\ref{mini}).
\end{enumerate}
\end{theo}

Moreover, it was shown that in the case where $c$ is the squared Riemannian distance on $M$, denoting $\Sigma_g$ the sectional curvature of $M$, the identity
\be
\mathfrak{S}_c(x,x)= \frac{3}{2} \Sigma_g(x)
\en holds (meaning that the two tensors coincide on $(T_x\Omega)^2$). Gathering all those results implies that for a manifold whose sectional curvature is negative at some point on some two-plan, one can exhibit positive smooth densities such that the optimal transport map between them (for the  cost equal to $d^2(\cdot,\cdot)$) is not continuous.

\section{Connectedness of the contact set: Proof of Theorem \ref{main_sphere_1}}

We prove the following more general result:
\begin{theo}\label{main_general}
Let $\Snn$ be the unit sphere of $\Rd$ equipped with the round metric $g$, and Riemannian distance $d$. Let $c(x,y)=f(d(x,y))$ for some $f: [0,\pi) \to \R$ smooth and strictly increasing, with $f'(0)=0$ and such that the conditions {\bf A0}, {\bf A1}, {\bf A2}, {\bf As} are satisfied on $\Snn\times \Snn \setminus \text{antidiag}$. Then 
\begin{enumerate}
\item for all $x, y_0, y_1 \in \Snn$ such that $-x \notin \{y_0, y_1\}$,  for all $y \in [y_0, y_1]_x$, $y\neq -x$;
\item for all $\phi$ c-convex, for all $y\in \Snn$, $\phi + c(\cdot, y)$ reaches a global maximum at $x=-y$, and any other critical point is a global minimum;
\item for all $\phi$ c-convex, for all $x\in \Snn$, $\partial\phi(x) = \partial^c\phi(x).$
\end{enumerate}
\end{theo}

\textsc{Remark.} We notice that $d(x,y) = \arccos(x\cdot y)$, and that
\be
-\log|x+y| &=& -\demi(\log(\demi|x+y|^2) + \text{cste}\\
&=&-\demi\log(1 + \cos(d(x,y))) + \text{cste}.
\en
Hence in both cases the cost-function is of the form $c(x,y)=f(d(x,y))$ that fits into the framework of Theorem \ref{main_general}.

\noindent
\textsc{Proof of Theorem \ref{main_general}.} 
We first prove the point 1: as the cost $c$ is under the form $c=f(d)$, we have, for all $x\neq y$, 
\beq\label{gradsphere}
\nabla_x c(x,y) = -f'(d(x,y))e_y,
\enq
where $e_y = \exp^{-1}_x(y)/d(x,y)$. If the cost satisfies {\bf A1}, then $f'$ should obviously be a strictly monotone function. This yields that $f$ is necessarily strictly convex. Now since $f$ is non-negative and increasing, from the definition of the c-segment, we  have easily using (\ref{gradsphere})
\be
\max\{d(x,y_\theta), \theta\in [0,1]\}\leq \max\{d(x,y_0), d(x,y_1)\}.
\en
This yields point 1. Note that point 1 implies that the c-segment $[y_0,y_1]_x$ is well defined provided $-x\notin\{y_0, y_1\}$.
We now prove point 2 and 3 for $\phi$ of the form 
\beq\label{simple}
\bar\phi(x)=\max\{-c(x,y_0) + a_0, -c(x,y_1) + a_1\}.
\enq
We let 
\be
h(x) &=& \bar\phi(x)+c(x,y),\\
D_y &=& \Snn \setminus \{-y\}.
\en
Note that $h$ is semi-convex on $D_y$. We say that $x\in D_y$ is a critical point of $h$ whenever $0\in \partial h(x)$. 
We now use the fact that the cost satisfies assumption {\bf As}: 
\begin{lemme}\label{locmin} Under the assumptions of Theorem \ref{main_general},  let $x\neq -y$ be a critical point of $h$, then  
$x$ is a local minimum of $h$. Moreover $h$ reaches its global maximum at $x=-y$. 
\end{lemme}
\textsc{Proof.} We choose $x\neq -y$ and  assume that $\bar\phi$ is not differentiable at $x$ otherwise the conclusion is trivial. Without loss of generality, one can then write 
\be
\bar\phi(x')=\max\{-c(x',y_0)+c(x,y_0), -c(x',y_1)+c(x,y_1)\}.
\en
Moreover, if $\bar\phi + c(\cdot, y)$ has a critical point at $x$ where $c(\cdot, y)$ is differentiable, then necessarily $y\in[y_0,y_1]_x$, the  c-segment with respect to $x$. Hence, one can write $y=y_\theta$ for some $\theta\in[0,1]$. Finally, if $y_0 \neq y_1$, $\bar\phi$ is necessarily differentiable at $x=-y_i, i=0,1$, since in a neighborhood of $y_0$ (resp. $y_1$), $\bar\phi \equiv -c(x',y_1)+c(x,y_1)$ (resp. $\bar\phi \equiv -c(x',y_0)+c(x,y_0)$). Hence we can assume that $x\neq -y_0$ and $x\neq -y_1$.
If $y=y_0$ or $y=y_1$ then  $\bar\phi + c(\cdot, y)$ has a global minimum at $x$. In the case where $y\notin \{y_0, y_1\}$
we use the result of \cite[Proposition 5.1]{L5}. 
Locally around $x$,  under assumption {\bf As}, we have for $\theta \in [\epsilon, 1-\epsilon]$  
\be
&&\max\{-c(x',y_0)+c(x,y_0), -c(x',y_1)+c(x,y_1)\}\\ 
&\geq& -c(x', y_\theta) + c(x,y_\theta) + \delta d^2(x,x') + o(d^2(x,x')),
\en
where $\delta>0$ for $\epsilon>0, d(y_1,y_0)>0$. Hence, we have
\be
\bar\phi(x')+c(x',y_\theta) \geq c(x, y_\theta) + \delta d^2(x,x') + o(d^2(x,x')),
\en
with equality at $x'=x$. So if $\bar\phi + c(\cdot, y)$ has a critical point at $x$, either $y=y_0$ or $y_1$ in which case $x$ is a global minimum, or $y=y_\theta$ for some $\theta \in (0,1)$, and then $x$ is a strict local minimum.

This implies that $h$ reaches its global maximum at $x=-y$: Indeed, if $h$ reaches its maximum at $x\in D_y$, then $h$ is semi-convex around $x$ (since it is the sum of a c-convex function and a smooth function), and $x$ is thus a critical point of $h$. As we saw in Lemma \ref{locmin}, any  critical point of $h$ in $D_y$ is a local minimum, and as we are on a compact manifold, $h$ must have a global maximum, which is thus necessarily $-y$.

 $\hfill\Box$
 
\begin{lemme}\label{onemin}
Let $h, \bar\phi$ be defined as above. Then any local minimum of $h$ is a global minimum on $\Snn$.
\end{lemme}
\textsc{Proof.} We first assume that $y\notin \{y_0, y_1\}$. In this case, we prove the stronger assertion that $h$ has only one local minimum (which in particular is the global minimum). So let us assume that $h$ has two distinct local minima $x_1\neq x_2$. Note that from Lemma \ref{locmin}, both $x_1$ and $x_2$ are different from $-y$, which is the global maximum of $h$. From the proof of Lemma \ref{locmin}, these two local minima are strict, and are thus separated: $x_1$ and $x_2$ do not belong to the same connected component of a level set of $h$. Then consider
\be
\inf_{\gamma\in\Gamma}\sup_{t\in [0,1]} {h(\gamma(t))},
\en
where $\Gamma$ is the set of continuous paths from $[0,1]$ to $\Snn$ such that $\gamma(0)=x_1, \gamma_1 = x_2$. By standard compactness arguments, this saddle point is attained at some point $x_3$. Clearly $x_3\neq -y$ since, from Lemma \ref{locmin}, $h$ reaches its global maximum at $-y$. Thus $x_3\in D_y$, and $x_3$ is a critical point of $h$. In view of Lemma \ref{locmin}, $x_3$ should then be a local minimum of $h$,  moreover, since $y\notin \{y_0, y_1\}$, this local minimum should be strict. We reach a contradiction.

If $y=y_0$, the situation is simpler:   $h$ reaches its global minimum on the whole set 
\be
E=\{x': -c(x',y_0)+ a_0\geq -c(x',y_1)+a_1\}.
\en
Moreover $h$ reaches its global maximum at $x'=-y_1$, and, from assumption {\bf A1}, there are no critical points outside of $E$, so there can be no local minima of $h$ outside of $E$.

$\hfill\Box$

\begin{lemme}\label{special} 
Let $\bar\phi$ be defined as above. Then for all $x\in \Snn$, $\partial^c\bar\phi(x) = \partial\bar\phi(x)$.
\end{lemme}
\textsc{Proof.} Consider  $p\in \partial\bar\phi(x)$. Then $p\in[-\nabla_xc(x,y_0),-\nabla_xc(x,y_1)]_x$. If $x$ is on the set where $\phi$ is not differentiable, (otherwise there is nothing to prove) then $-x \notin \{y_0, y_1\}$. From the point 1 already proved, there exists $y\in [y_0, y_1]_x \subset \Snn\setminus \{-x\}$, such that $p=-\nabla_xc(x,y)$. Then considering the function $h = \bar\phi + c(\cdot,y)$ on $\Snn$, clearly $x$ is a critical point of $h$. Thus, by what we proved above, $x$ is a global minimum of $h$, and so $p\in \partial^c\bar\phi(x)$. 

$\hfill\Box$

Then arguing as in \cite[Proposition 2.11]{L5}, if the conclusion of Lemma \ref{special} holds for all c-convex functions of the form (\ref{simple}), it holds for all $\phi$ c-convex; this achieves the proof of Theorem \ref{main_general}.

$\hfill\Box$

\textsc{Remark.} This argument allows to go from  \cite[Proposition 5.1]{L5} to the important identity $\partial^c\phi = \partial\phi$. One of the crucial ingredients is the absence of boundary used in the saddle point argument. The case with boundary has been treated by Trudinger and Wang in \cite{TW2}. An alternative argument to obtain  $\partial^c\phi = \partial\phi$ from {\bf Aw} has been obtained in \cite{KimMcCann}.

\section{Proof of Theorem \ref{main_sphere}}

\paragraph{Strategy of the proof}
Most of the proof is contained in the following points:

1- Given $\mu_0, \mu_1$ satisfying the assumptions of Theorem \ref{main_sphere}, there exists a constant $\sigma$ such that $d(x,G_\phi(x)) \leq \pi -\sigma$ for all $x \in \Sn^{n-1}$. Hence, $G_\phi(x)$ stays uniformly far away from the cut-locus of $x$. Then we can reduce locally the problem to an Euclidean problem.

2- The assumption {\bf As} is satisfied by the cost function distance squared on the sphere (Theorem \ref{Assphere}).
%\begin{theo}
%Let $\Sn^{n-1}$ be the unit sphere of $\Rd$ equipped with the round metric $g$, and Riemannian distance $d$. Let  $c(x,y)=\demi d^2(x,y)$. Then $c$ satisfies {\bf As} on  $\Sn^{n-1}\times \Sn^{n-1}\setminus \{(x,x), x\in \Sn^{n-1}\}$. In other words, for this choice of $c$, the cost-sectional curvature of $\Sn^{n-1}$ is uniformly positive.
%\end{theo}

Once this is established, we proceed as follows:

Given $x_0 \in \Sn^{n-1}$ we can build around $x_0$ a system of geodesic coordinates on the set $\{x, d(x,x_0) \leq R\}$ for $R<\pi$. From point 1-, for $r$ small enough, the graph $\{(x, G_\phi(x)), x\in B_r(x_0)\}$ is included in the set $B_r(x_0)\times B_{\pi-2r}(x_0)$ on which the cost function is $C^\infty$. From point 2- and using \cite{MTW}, a $C^4$ smooth solution to (\ref{ma}) on $B_r(x_0)$ will enjoy a $C^2$ a priori estimate at $x_0$. This estimate will depend only  on the smoothness of $\mu_0, \mu_1$, on $r$, and $r$ is small but can be chosen once for all.

Once a $C^2$ a priori estimate is established we use the  result of \cite{del2}:  the method of continuity allows to build smooth solutions for any smooth positive densities. 

Then, we use the results of \cite{L5} to conclude our partial $C^{1,\alpha}$ regularity result.

%As in \cite{L5}, we can infer from the existence of smooth solutions that $\partial\phi = \partial^c\phi$.

\subsection{Reduction of the problem to an Euclidean problem} 

%We denote by $B_r$ (resp. $B_r(x)$) a Riemannian ball of radius $r$ (resp. centered at $x$), $\Sn^{n-1}$ denotes the unit sphere of $\Rd$.

\paragraph{Uniform distance to the cut locus}
We show that there exists a subset $S^2_\sigma$ of $\Sn^{n-1} \times \Sn^{n-1}$ on which {\bf A0}-{\bf A2} are satisfied, and such that the graph of $G_\phi$, i.e. $\{(x, G_\phi(x)), x\in \Sn^{n-1}\}$, is contained in $S^2_\sigma$. This subset $S^2_\sigma$ is defined by 
\beq\label{defS2sigma}
S_\sigma^2=\Big\{(x,y) \in \Sn^{n-1} \times \Sn^{n-1}: d(x,y) \leq \pi - \sigma\Big\}
\enq  
where $\sigma >0$ depends on some condition on $\mu_0, \mu_1$. 
We first remark that 
\begin{lemme}\label{qdmmccvx}
For all $x\in \Sn^{n-1}$, the set $\{y, (x,y) \in S_\sigma^2\}$ is c-convex with respect to $x$.
\end {lemme}
Then we show the  following crucial result:
\begin{prop}\label{sphereok}
Let $\mu_0, \mu_1$ be two probability measures on $\Sn^{n-1}$, let $\phi$ be a c-convex potential such that $G_{\phi\#}\mu_0=\mu_1$. Assume that there exists $m>0$ such that $\mu_1 \geq m{\rm dVol}$ and that $\mu_0$ satisfies {\bf B(n-1)} (see (\ref{mini})). Then there exists $\sigma >0$ depending on $m$ and on   $f$ in (\ref{mini}), such that $\{(x, G_\phi(x)), x\in \Sn^{n-1}\} \subset S^2_\sigma$, where $S^2_\sigma$ is defined in (\ref{defS2sigma}).
\end{prop}

\textsc{Proof:} 
We use \cite{del-loep}; in that paper, it was shown, for $\phi$ satisfying $G_{\phi\#}\mu_0=\mu_1$, 
that we have $|d\phi| \leq \pi - \epsilon$, and $\epsilon>0$ depends on  $\|{\rm d}\mu_1/{\rm dVol}\|_{\Linf}\|\left({\rm d}\mu_0/{\rm dVol}\right)^{-1}\|_{\Linf}$.
Considering $\phi^c$ and $G_{\phi^c}$ that pushes forward $\mu_1$ onto $\mu_0$,  we also have a bound such as $|d\phi|\leq \pi -\epsilon$, where $\epsilon>0$ depends on 
$\|{\rm d}\mu_0/{\rm dVol}\|_{\Linf}\|\left({\rm d}\mu_1/{\rm dVol}\right)^{-1}\|_{\Linf}$. Here we slightly extend this bound to the case where $\mu_0$ satisfies {\bf B(n-1)} (\ref{mini}), and $\mu_1\geq m{\rm dVol}$.

The starting point of the proof is  \cite[Lemma 2]{del-loep} where it is shown that for all $x_1, x_2 \in \Sn^{n-1}$, for all $\psi$ c-convex,
\be
d(G_\psi(x_2), \hat x_1) \leq 2\pi\frac{d(G_\psi(x_1), \hat x_1)}{d(x_1,x_2)},
\en
where $\hat x_1$ is the antipodal point to $x_1$.
Hence the set $E_\delta=\{x: d(x,x_1)\geq \delta\}$ is sent by $G_\psi$ in $B_\epsilon(\hat x_1)$ where
$$\epsilon = 2\pi\frac{d(G_\psi(x_1), \hat x_1)}{\delta}.$$  
Considering $\psi = \phi^c$, we have $G_{\psi\,\#}\mu_1 = \mu_0$.
This implies
\be
\mu_0 (B_\epsilon(\hat x_1)) \geq \mu_1 (E_\delta).
\en
Taking $\delta > 0$ fixed (for example $\delta = \pi/2$), we have
\be
\mu_0(B_{4d(G_{\phi^c}(x_1), \hat x_1)}(x_1)) \geq \inf_z \mu_1(D_z),
\en
where $D_z$ is the half-sphere centered at $z$.
Hence
\be 
&&\inf_{x\in \Sn^{n-1}}\Big\{|\pi- d\phi^c(x)|\Big\}\\
&\geq& \inf\Big\{r: \exists x\in \Sn^{n-1}, \ \mu_0(B_{4r}(x)) \geq  \inf_z \mu_1(D_z) \Big\}\\
&:=&\sigma.
\en
In particular,  if $\mu_0$ satisfies {\bf B(n-1)},   we have  $$\mu_0(B_{4r}(x)) \leq  C\left({\rm Vol}(B_{4r}(x))\right)^{1-1/(n-1)}$$ for $r$ small enough,
hence if  $\mu_1 \geq m{\rm dVol}, m>0$ we find that $\sigma>0$  and the uniform distance 
between $G_\psi(x)$ and the cut locus  of $x$ is asserted. Now notice that $\sup\{|d \phi(x)|, x\in\Snn\} = \sup \{|d \phi^c(x)|, x\in\Snn\}$, and the  proposition follows.

$\hfill \Box$

\textsc{Remark.} It is enough to conclude the proof to assume that $\mu_0(B\epu(x))\to 0$ uniformly with respect to $\epsilon$.

\paragraph{Construction of a local system of coordinates}
Given $x_0\in \Sn^{n-1}$, we consider a  system of geodesic coordinates around $x_0$, i.e. given a system of orthonormal coordinates at $x_0$ and  the induced system of coordinates on $T_{x_0}(\Sn^{n-1})$, we consider the mapping $p\in T_{x_0}(\Sn^{n-1}) \to \exp_{x_0}(p)$. 
This mapping is a diffeomorphism from $B_R(0) \subset \R^{n-1}$ to $B_R(x_0)$ as long as $R<\pi$. Then, for for $r<\sigma/2$ we have 
$$S^2_\sigma \cap \{(x,y): x\in B_r(x_0)\} \subset B_r(x_0) \times B_{\pi-\sigma + r} \subset S^2_{\sigma - 2r}.$$
We now take $r=\sigma/3$
Hence,  $B_r(0)$ is sent in $B_{\pi-2r}(0)$, and the cost function is $C^\infty$ on $B_r(0) \times B_{\pi-2r}(0)$for $r$ small enough.

In this case, we have $\mT_x(p) = \exp_x(p)$, and also 
\beq\label{note}
D^2_{x,y} c(x,y) = -[D_v \exp_x]^{-1}(y),
\enq
where $D_v\exp_x$ is the derivative with respect to $v$ of $v\to \exp_x(v)$.
Assumption {\bf A1} is trivial, since $y$ is indeed uniquely defined by $y=\exp_x(p)$. From (\ref{note}), assumption {\bf A2} is true on any smooth compact  Riemannian manifold, since in this case  ${\rm Jac}(v \to \exp_x v)$ is bounded by above.

\subsection{Verification of assumption {\bf As}: Proof of Theorem \ref{Assphere}}\label{sphere}

It has been established in \cite{del-loep} that in a system of normal coordinates $e_1,..,e_{n-1}$ at $x$ with $e_1= p/|p|$, we have
\be
\ds D^2_{x_ix_j}c|_{x,\exp_x(p)} = \left\{\begin{array}{ll} 1  \text{ if } i=j=1 \\\ds \delta_{ij} \frac{|p| \cos(|p|)}{\sin(|p|)} \text{ otherwise}. \end{array}\right.
\en
This relies on the fact that $D^2_{xx}c(x,\exp_x(p))=[D_p(\exp_x(p))]^{-1}D_x(\exp_x(p))$, and follows by computations of Jacobi fields.
Hence, we have, in this system of coordinates,
\be
M(x,p):=D^2_{xx}c(x, \mT_x(p)) = I - (1-\frac{r\cos r}{\sin r})I',
\en
where $r=|p|$, $I$ is the identity matrix of order $n-1$, and $I'$ is the identity matrix on $\text{vect}(e_2,..,e_{n-1})$.
So, for a given $v \in \R^{n-1}$,
\be
M(x,p)\cdot v\cdot v = |v|^2 - (1-\frac{r\cos r}{\sin r})|\Pi v|^2,
\en
where $\Pi v$ is the projection of $v$ on $p^\perp$.
This can be written, using intrinsic notations 
\be
\left({\cal H}(c(\cdot, y))|_{y=\exp_{x}(z)}v,v\right)_g=
|v|^2_g - (1-\frac{|z|_g\cos |z|_g}{\sin |z|_g})|\Pi v|_g^2,
\en
where $(\cdot,\cdot)_g$ denotes the Riemannian scalar product on $\mT_x(\Sn^{n-1})$ and $|p|_g^2=(p,p)_g$.
 
Let $\alpha , \beta \in \R^{n-1}$, $t\in \R$, $u_t=\alpha+ t\beta$ and $r=|u_t|$. We assume that $r < \pi$, and $\Pi_t v$ is the projection of $v$ on
$u_t^{\perp}$.
Assumption {\bf As} is then equivalent to show
\be
\Dtt \left[t\to M(x,\alpha+\beta t)\cdot v\cdot v \right]\leq -C_0 |\beta|^2 |v|^2,
\en
for all $\alpha, \beta \in \mT_x\Sn^{n-1}$, $v\perp \beta$, which is equivalent to
\be
\Dtt \left[t \to (1-\frac{r\cos r}{\sin r})  |\Pi_t v|^2\right] \geq C_0 |\beta|^2 |v|^2,
\en
for all $\alpha, \beta \in \R^{n-1}$, $v\perp \beta$. Without loss of generality, we assume by changing $t$ in $t-t_0$ that $\alpha \perp \beta$.
We have
\be
\Pi_t v = v - |u_t|^{-2} (v\cdot u_t) u_t,
\en
hence, using $v\perp\beta$,
\be
|\Pi_t v|^2 &=& |v|^2 - 2 (v\cdot u_t)^2|u_t|^{-2} + (v\cdot u_t)^2|u_t|^{-2}\\
&=& |v|^2 - (v\cdot \alpha)^2|u_t|^{-2}.
\en
Now, assuming $|v|=1$ and $ (v,\alpha)_g =c$, we have to evaluate $\Dtt F$, where $F=G\circ r$,
and
\be
&&G(r) = (1-\frac{r\cos r}{\sin r})(1- \frac{c^2}{r^2}),\\
&&r(t) = |\alpha+ t\beta|.
\en
A computation shows that
\be
\Dtt F &=&\frac{\alpha^2 \beta^2}{r^3} \Bigg\{\frac{1}{\sin^2 r}(r-\sin r\cos r)(1 - \frac{c^2}{r^2})      \\ 
&& \hspace{1cm} + \frac{1}{\sin r}(\sin r - r\cos r)\frac{2 c^2}{r^3} \Bigg\}\\
&+& \frac{t^2\beta^4}{r^2} \Bigg\{\frac{2}{\sin^3 r}(\sin r -r\cos r)(1-\frac{c^2}{r^2}) \\
&&\hspace{1cm} + \frac{2}{\sin^2 r}(r - \sin r \cos r ) \frac{2 c^2}{r^3} \\
&&\hspace{1cm} - \frac{1}{\sin r}(\sin r -r \cos r) \frac{6 c^2}{r^4} \Bigg\}.
\en
Notice that we have for some $C>0$, and for all $r\in [0,\pi]$, 
\beq
&&\sin r - r\cos r \geq Cr^3,\label{util1}\\
&&r-\sin r\cos r \geq Cr^3\label{util2}.
\enq
Hence,  all terms are positive, except the last one.
The sum of the last two lines has for $r\in[0,\pi]$ the sign of 
\be
&&2(r^2-r\sin r\cos r) - 3(\sin^2 r -r\sin r \cos r)\\
&=&\demi \left(r\sin 2r + 3 \cos 2r + 4 r^2 -3 \right).
\en
One can check (by two successive differentiations) that this quantity  is non-negative for $r\in [0,\pi]$.
Hence we have, using (\ref{util1}, \ref{util2}), 
\be
\Dtt F  &\geq& C\beta^2 \left[(1 - \frac{c^2}{r^2})  \frac{\alpha^2 + t^2\beta^2}{r^2} + \frac{\alpha^2 c^2}{r^4}\right]\\
&=& C\beta^2\left[(1-\frac{c^2}{r^2})+ \frac{\alpha^2 c^2}{r^4}\right].
\en
Remember that $c = \alpha\cdot v$, with $|v|=1$, hence  $|\alpha| \geq |c|$ and $r\geq |c|$, to conclude that $\Dtt F  \geq C_0|\beta|^2$.

Hence in a normal system of coordinates at $x$, we see that {\bf As} is satisfied at $x$ for any $y$ such that $d(x,y) < \pi$. As seen in \cite{L5}, the cost-sectional curvature is intrinsic. This proves Theorem \ref{Assphere}.

$\hfill \Box$

\subsection{$C^1$ regularity for weak solutions}
%This follows by following the proof of Theorem \ref{main} in \cite{L5}. 
Consider a system of geodesic coordinates  around $x_0$.  
As we have seen, using Proposition \ref{sphereok},  in this system of coordinates, for $r_0=\sigma/3$, we have $$\{(x, G_\phi(x)), x\in B_{r_0}(0)\} \subset B_{r_0}(0) \times B_{\pi-2r_0}(0) \subset S^2_{\sigma/3}.$$ 
%In this system of coordinates,  $\phi$ will satisfy on $B_{r_0}(0)$ the Monge-Amp\`ere equation (\ref{ma}) where $c$ is $d^2/2$ expressed in our local coordinates. Moreover, 
Moreover,  $c$ satisfies {\bf A0}-{\bf As} on $B_{r_0}(0) \times B_{\pi-2r_0}(0)$.
We let $\Omega = B_{r_0}(0)$, $R=\pi-2r_0$, and $\Omega'=B_R(0)$, and remark that $\Omega, \Omega'$ are uniformly strictly c-convex with respect to each other if $r_0$ is small enough.
Now, given $\phi$ c-convex on $\Snn$ such that $G_{\phi\,\#}\mu_0=\mu_1$, we consider $\tilde\phi $ the restriction of $\phi$ to $\Omega$. 
We recall the definition of the Monge-Amp\`ere measure of $\phi$ with respect to $\mu_1$, defined by
\beq
\forall B \text{ Borel}, G_\phi^{\#}\mu_1(B) = \mu_1(G_\phi(B)).
\enq

We claim the following:
\begin{prop}\label{restriction}
Let  $\phi, \tilde\phi, \Omega, \Omega'$ be as above. Let $\mu_0, \mu_1$ be probability measures on $\Snn$ such that $G_{\phi\,\#}\mu_0 = \mu_1$, and $\mu_1\geq m \rm{dVol}$; then 
\begin{enumerate}
\item $\tilde\phi$ is c-convex on $\Omega$ with respect to the restriction of $c$ to $\Omega\times\Omega'$;
\item $G_{\tilde\phi}$ coincides with $G_\phi$ on $\Omega$;
\item $G_{\tilde\phi}^{\#}\rm{dVol}\leq \frac{1}{m}\rm{dVol}$ on $\Omega$.
\end{enumerate}
\end{prop}
\textsc{Remark.} Note that in general, $G_{\phi\,\#}\mu_0{\bf 1_\Omega} =  \mu_1{\bf 1_{\Omega'}}$ will not hold. Indeed, $G_\phi(\Omega)$ is a priori strictly included $\Omega'$,  hence for $B\subset\Omega'$, $\mu_1(B)=\mu_0(G_\phi^{-1}(B))$ will be strictly larger than $\mu_0(G_\phi^{-1}(B)\cap \Omega)$. However, point 2. suffices to conclude that $G_{\tilde\phi}^{\#}\mu_1 = G_{\phi}^{\#}\mu_1$.

\textsc{Proof.} The first point is straightforward: indeed at every $x\in\Omega$, $\tilde\phi$ has global c-support $-\phi^c(y) - c(\cdot,y)$, and $y\in \Omega'$ since $G_\phi(\Omega)\subset\Omega'$.
  
The delicate issue for the second part comes from the following point: we have
\be
G_{\tilde\phi}(x)=\Big\{y: \phi(x) + c(x,y)= \inf\{\phi(x') + c(x',y), x'\in \Omega \}\Big\}.
\en
As we see, the infimum above has to be on $\Omega$, whereas for the definition of $G_\phi$ it has to be on $\Snn$.
Thus the inclusion $G_\phi(x) \subset G_{\tilde\phi}(x)$ is straightforward while the reverse inclusion is unclear.
To overcome this we use the Theorem \ref{main_general}, and the fact that for all $\phi$ c-convex on $\Snn$, $\partial^c\phi = \partial\phi$. Hence if $\phi + c(\cdot,y)$ has a local minimum  it has to be a global minimum in all of $\Snn$.
Hence, if $y\in G_{\tilde\phi}(x)$, then $\phi + c(\cdot,y)$ has a local minimum at $x$, and hence a global one, thus $y\in G_\phi(x)$.

Now the last point follows from the following lemma, which is  a reformulation of \cite[Proposition 5.11]{L5} in the sphere case.

\begin{lemme}\label{##}
Let $\phi$ be c-convex on $\Snn$. Assume that $G_{\phi\,\#}\mu_0=\mu_1$. Assume that $\mu_1\geq m{\rm dVol}$ on $\Omega'$. Then for all $\omega \subset \Snn$, we have
\be
\mu_0(\omega) \geq m {\rm Vol}(G_\phi(\omega)),\text{ and hence, } G_\phi^{\#}{\rm dVol} \leq \frac{1}{m} \mu_0.
\en
\end{lemme}

Before proving this lemma, we conclude the proof of Proposition \ref{restriction}, by observing that, from the point 2 of Proposition \ref{restriction},  on $\Omega$ we have $G_{\tilde\phi}^{\#}{\rm dVol}= G_\phi^{\#}{\rm dVol}$.

$\hfill\Box$

\textsc{Proof of Lemma \ref{##}.} In $\Snn$ we consider $N=\{y\in \Snn: \exists x_1\neq x_2 \in \Snn, G_\phi(x_1)=G_\phi(x_2)=y\}$. Then $N=\{y\in \Snn: \phi^c \text{ is not differentiable at } y\}$. Hence ${\rm Vol}(N)=0$, and  ${\rm Vol} (G_\phi(\omega)\setminus N) = {\rm Vol}(G_\phi(\omega))$. Moreover,  on $G_\phi(\omega)\setminus N$, $G_\phi^{-1}$ is single valued. Then  $G_\phi^{-1}(G_\phi(\omega)\setminus N) \subset \omega$. Hence, 
\be
\mu_0(\omega) &\geq& \mu_0 (G_\phi^{-1}(G_\phi(\omega)\setminus N))\\
&=&\mu_1(G_\phi(\omega)\setminus N)\\
&\geq & m{\rm Vol}(G_\phi(\omega)\setminus N)\\
& = & m{\rm Vol}(G_\phi(\omega)).
\en

We are now allowed to use \cite[Proposition 5.9]{L5} to prove Theorem \ref{main_sphere}. Let first recall the result:

\begin{prop}[Proposition 5.9 \cite{L5}] Let $\Omega, \Omega'$ be domains of $\R^{n-1}$. Let $\phi$ be c-convex with $G_\phi(\Omega)\subset \Omega'$. Let $c$ satisfy assumptions  {\bf A0}, {\bf  A1}, {\bf A2}, {\bf As} on $\Omega, \Omega'$ with $\Omega, \Omega'$ strictly c-convex with respect to each other, then
 \begin{itemize}
\item 
if $G_\phi^{\#}{\rm dVol}$ satisfies {\bf A(n-1, p)} (see (\ref{Lp})), for some $p>n-1$, then $\phi \in C^{1, \beta}_{loc}(\Omega)$, with $\beta(n-1,p)$ as in Theorem \ref{main},
\item
if $G_\phi^{\#}{\rm dVol}$ satisfies {\bf B(n-1)} (see(\ref{mini})), then $\phi \in C^1_{loc}(\Omega)$.
\end{itemize}
%where $G_\phi^\#$ is defined in Definition \ref{defiweak}.

\end{prop}
This achieves the proof of $C^1$ regularity.

\subsection{$C^2$ a priori estimates for smooth solutions}

Let $\phi$ be a $C^4$ smooth c-convex potential on $\Snn$ such that $G_{\phi \, \#}\mu_0 = \mu_1$, with $\mu_0, \mu_1$ having positive $C^{1,1}$ smooth densities. We recall the a-priori estimate of  \cite[Theorem 4.1]{MTW}, adapted to our notations:
\begin{theo}[\cite{MTW}] Let $\Omega, \Omega'$ be subsets of $\Rd$. Let $\phi$ be a $C^4(\Omega)$ c-convex solution of (\ref{ma}) with $G_\phi(\Omega)\subset\Omega'$. Suppose assumptions
{\bf A0}, {\bf A1}, {\bf A2}, {\bf As} are satisfied by $c$ on $\Omega\times\Omega'$. Then we have the a priori second order derivative estimate
\be
|D^2\phi|(x)\leq C,
\en
where C depends on n, $dist(x, \partial\Omega)$, $\sup_{\Omega}\phi$, $\rho_0, \rho_1$ up to their second derivatives, the cost
function c up to its fourth order derivatives, the constant $C_0$ in {\bf As}, and a positive
lower bound of $|\det D^2_{x,y}c|$ on $\Omega$.
\end{theo}
Applying this result with $\Omega = B_r(0)$, $\Omega' = B_{\pi-2r_0}(0)$, we have a bound on $D^2\phi(x_0)$ that depends only on $r_0$ and on the bounds on $\mu_0, \mu_1$. Finally, note  that $r_0$ can be chosen once for all once $\sigma$ is known.

\paragraph*{Continuity method}
In \cite{del2} it was established that given a $C^\infty$ smooth  c-convex potential $\phi$, a $C^\infty$ smooth positive measure $\mu_0$, and $\mu_1=G_{\phi\,\#}\mu_0$, the operator
\be
F: \phi\to F(\phi)=G_{\phi\,\#}\mu_0
\en 
was locally invertible in $C^\infty$ around $\mu_1$. 
Then the existence, for a given pair of $C^\infty$ smooth positive probability measures $\mu_0,\mu_1$ of a $C^\infty$ smooth c-convex potential $\phi$ is granted once a-priori estimates for the second derivatives have been derived. Indeed, this follows from the concavity of the equation and and the well known continuity method (see \cite{GT}).

We conclude the following: for $\mu_0, \mu_1$ having $C^{1,1}$ smooth probability densities, there exists a (unique up to a constant) c-convex potential $\phi\in C^{3,\alpha}$ for every $\alpha \in [0,1[$ such that $G_{\phi\,\#}\mu_0=\mu_1$.
If moreover $\mu_0, \mu_1$ are $C^\infty$, $\phi \in C^\infty$.
This concludes the third  point of Theorem \ref{main_sphere}.

%\paragraph*{Locally supporting functions are globally supporting functions}
%Again, one we have proved existence of smooth solutions for arbitrary smooth positive data, this follows from \cite{L5}. More precisely, we have
%\begin{theo}
%For any $\phi$ c-convex on $\Sn^{n-1}$ with $c = d^2/2$, for all $x \in \Sn$, $\partial\phi(x) = \partial^c\phi(x)$. 
%\end{theo}

%Alternatively, once the existence of smooth solutions is established, the  theorem \ref{main} can be used as a $C^{1,\alpha}$ a priori estimate when approximating the data by smooth densities. As the bounds in Theorem \ref{main} will be uniform, one can easily pass to the limit. 
$\hfill \Box$

The proof of Theorem \ref{main_sphere} is complete.

\section{Proof of Theorem \ref{mainantenna}} For this case it has already been checked (see \cite{Wa3, Wa4, GuWa})  that the cost function $c(x,y)=-\log|x-y|$ satisfies {\bf As} for $x\neq y$. We will just prove that under some assumptions on the measures $\mu_0, \mu_1$ we can  guarantee that $G_\phi(x)$ stays away from $x$. 
This will imply that when the cost is $-\log|x+y|$, $G_\phi(x)$ stays away from $\hat x$.
Then, the proof of the Theorem \ref{mainantenna} mimics the proof of Theorem \ref{main_sphere}.
We prove the 
\begin{prop}\label{awayantenna}
Let $\Sn^{n-1}$ be unit sphere of $\Rd$. Let $c(x,y) = -\log|x-y|$.
Let $T:\Sn^{n-1}\to \Sn^{n-1}$ be a 2-monotone map, i.e. such that
\beq\label{antenna2monot}
\forall x_1, x_2 \in \Sn^{n-1}, c(x_1, T(x_1)) + c(x_2, T(x_2)) \leq c(x_2, T(x_1)) + c(x_1, T(x_2)).
\enq
Let $\mu_0, \mu_1$ be two probability measures on $\Rd$. Assume that $\mu_1\geq m{\rm dVol}$ for $m>0$, and that $\mu_0$ satisfies {\bf B(n-1)} (\ref{mini}). 
Then there exists $\epsilon_0 >0$ depending only on $m$, $f$ in (\ref{mini}) such that 
\be
\forall x\in \Sn^{n-1}, d(x,T(x)) \geq \epsilon_0.
\en
\end{prop}

\textsc{Proof of Proposition \ref{awayantenna}.} We follow the same lines as in \cite{del-loep}.
From (\ref{antenna2monot}), we have
\be
\log|x_1-T(x_2)| \leq \log|x_1-T(x_1)| + \log 2 - \log|x_2-T(x_1)|.
\en
Letting $M=-\log|x_1-T(x_1)|$, the set $\{x: |x-T(x_1)| \geq 1\}$ is sent by $T$ in the set $\{y: |y- x_1| \leq 2\exp(-M)\}$. Then the bounds on $\mu_0, \mu_1$ yield an upper bound on $M$ as in the proof of Proposition \ref{sphereok}.

$\hfill \Box$

\textsc{Remark.} The results that we have established for the reflector antenna and the quadratic cost can certainly be extended to the class of all costs that follow the assumptions of Theorem \ref{main_general}. Indeed the argument used to show the "stay away from the cut locus" property can be easily adapted in this case, and all the other arguments remain true.

\bibliography{C1-biblio}
\vspace{1cm}
\begin{flushright} 
Gr\'egoire Loeper
\\ 
Laboratoire Jacques-Louis Lions
\\
Université Pierre et Marie Curie
\\
Boîte courrier 187
\\
75252 Paris Cedex 05
\\  
France 
\\
loeper@ann.jussieu.fr
\end{flushright}

\end{document}